\documentclass{amsart}

\theoremstyle{plain}
\newtheorem{theorem}{Theorem}

\newtheorem{lemma}{Lemma}

\theoremstyle{remark}

\theoremstyle{definition}

\usepackage{multicol}
\usepackage{enumitem}

 \newcommand{\dd}[2]
{
{\frac{\partial #1}{\partial #2}}
}

\begin{document}

\title{Sard's conjecture for degenerate Engel}
\author{A. Perico}
\email{aperico@ucsc.edu}
\address{Department of Mathematics, University of California, Santa Cruz}


\date{September 2017}

\begin{abstract}
	Let $D$ be a rank $2$ bracket generating distribution on a $4$ manifold, $D$ is Engel if its growth vector is maximal. When this maximality fails the distribution is degenerate. We prove Sard's conjecture for the endpoint map in the case of degenerate Engel distributions. In this case the set of singular horizontal curves starting from the same point has measure zero: a full 2-dimensional disk.
\end{abstract}

\maketitle

Sard's conjecture for the endpoint map is an open problem in subRiemannian geometry. Here we solve the conjecture for a special class of cases: the generic degenerations of an Engel distribution. Previous work in the same spirit was done by Zelenko and Zhitomirskii \cite{zelenko1995rigid} for the generic degeneration of the Martinet distribution.

Let us recall that a subRiemannian geometry on a manifold $Q$ consists of a distribution (subbundle of the tangent bundle) $D\subset TQ$, together with a inner product on that distribution.

Let $\Omega$ be the space of absolutely continuous and square integrable paths $\gamma : [0, 1] \to Q$ such that $\gamma'(t)\in D_{\gamma(t)}$. For $q_0 \in Q $ consider $\Omega(q_0) = \{\gamma\in\Omega | \gamma(0) = q_0\}$. The endpoint map based at $q_0$ is the map 

\begin{align*}
	 end : \Omega (q_0) &\to  Q \\
	  		end (\gamma) &= \gamma(1)
\end{align*}

$\Omega(q_0)$ has the structure of a Hilbert manifold modeled on $L^2(I,\mathbb{R}^k)$ where $k$ is the rank of $D$ [Appendix \cite{montgomery1995survey}]. The map $end$ is smooth from an infinite dimensional manifold to the finite dimensional manifold $Q$. The Sard conjecture asserts that the Sard theorem holds for $end$, provided we restrict $end$ to short enough curves, specifically exists $C>0$ (dependent on the inner product) with the following property. Let $\Omega^C$ denote the open set of paths in $\Omega (q_0)$ of length less or equal than $C$, where the length is given by $l(\gamma)=\int_{\gamma}\sqrt{<\dot{\gamma},\dot{\gamma}>}$, the conjecture asserts that Sard holds for the restriction of $end$ to $\Omega^C$ it remains an open conjecture.

In more detail we also know from Bryan-Hsu \cite{bryantHsu1993rigidity} that a curve $\gamma$ is singular for $end$ if and only if there exists a covector $\lambda$ along $\gamma$ such that

\begin{enumerate}[label=\textbf{S.\arabic*}]
	\item \label{s1} $\lambda(t)\neq 0$ (almost everywhere).
	\item \label{s2} $\lambda(t)\in D^\perp$, the annihilator of $D$.
	\item \label{s3} $\lambda(t)$ is a characteristic curve: $\omega_{D^\perp}(\lambda(t))(\dot{\lambda}(t),\cdot)=0$; where $\omega_{D^\perp}$ is the restriction to the smooth subbundle $D^\perp\subset T^*Q$ of the canonical 2-form on $T^*Q$.
\end{enumerate}

The Sard conjecture asserts that the endpoints $\gamma(1)$, as $\gamma$ varies over the singular curves in $\Omega(q_0)$ has measure zero.

In order to describe the generic degeneration of an Engel distribution we need to recall the notion of the growth vector of a distribution at a point. We focus attention on rank two distributions $D$. $D$ is locally spanned by two vector fields $X, Y$. We can compute their brackets $[X,Y], [X,[X,Y]], \ldots$ Then, by definition $D^2$ is spanned by $X(p), Y(p)$ and $[X,Y](p)$ while $D^3$ is spanned by $X(p)$, $Y(p)$, $[X,Y](p)$, $[X,[X,Y]](p)$ and $[Y,[X,Y]](p)$. We continue in this manner, with $D^{k+1}$ spanned by $K$-fold bracketing $X$ and $Y$. The list of dimensions $dim D_p =n_1=2$, $dim D^2_p=n_2\geq 2$, $dim D^3_p=n_3$, \ldots, is called the growth vector of $D$ at $p$. If the manifold $Q$ is four dimensional and $D$ is bracket generating then the growth vector has the form $(2,n_2,n_3, \ldots , n_{s-1}, 4)$ with $n_i= 2$ or $3$ and $n_i\leq n_{i+1}$ for $i=2,3,\ldots, s-1$.

The condition of being Engel is that the growth vector is $(2,3,4)$ at every point. In this case there is a local normal form
\[ dz -xdw = 0,\quad dy-ydw=0\]
for the distribution.
The Engel distribution has a distinguished line field $L$, characterized by the condition
\begin{equation}
[L,D^2]\subset D^2.
\label{linefield}
\end{equation}
in the standard case above $L =span\{ \dd{}{w} \}$.

The integral curves  of $\dd{}{w}$ are exactly the singular curves of an Engel distribution.

To degenerate an Engel distribution we must violate the growth vector condition. The two generic degenerations are $(2,2,4)$ and $(2,3,3,4)$. Both occur along surfaces $\Sigma$ of $Q$. Off $\Sigma$, the distribution is Engel and so has a characteristic direction, through each point.

\begin{theorem}
	Let $q_0$ lie in $\Sigma$, one of the two surfaces of generic degeneration of an Engel distribution. Then Sard theorem holds at $q_0$. Indeed the set of critical values form a smooth $2$-disc   through $q_0$ in cases $(2,2,4)$ and $(2,3,3,4)$ case (A) below,
	and consists of the single point $\{q_0\}$ in the remaining case.  
	\end{theorem}

Zhitsomirskii \cite{zhitomirskii1990normal} found the normal forms for the degenerations of Engel, 
\begin{align*}
	(2,2,4) \text{ case:}& &\theta_1 &=dx + z^2dw, & \theta_2 &=dy + zwdw \\
	(2,3,3,4) \text{ case (A):}& &\theta_1&=dx+zdw, & \theta_2&= dy+ z^2wdw\\
	(2,3,3,4) \text{ case (B):}& &\theta_1&=dx +zdw, & \theta_2&= dy+ \left( \frac{1}{3} z^3+ zw^2 \right) dw 
\end{align*}

In these normal form coordinates, $q_0$ is the origin. The surface $\Sigma$ is given by $z=w=0$.
Off   $\Sigma$ the distribution $D$ is Engel.

The strategy of proof is to compute the characteristic direction $X$ (the rectification of the line field (\ref{linefield})) for our normal forms off of $\Sigma$ and show that the set of flow lines which end in the origin has dimension $2$. 

\begin{proof}
To start off, for any rank two distribution $D$ in a four manifold, we can find local coordinates $x, y, z, w$ and smooth functions $f$ and $g$ such that $D$ can be locally expressed as the Pfaffian system :

\begin{equation}
\tag{Pf}
\begin{split}
\theta_1&=dx +fdw\\
\theta_2&=dy +gdw
\end{split}
\label{eqn:Pfaffian}
\end{equation}

Using (\ref{eqn:Pfaffian}), we compute that 
$$D= span\{Z, W\} \text{ where }  Z = \dd{}{z}, W= \dd{}{w}-f\dd{}{x}-g\dd{}{y}.$$
We compute that $D$ is Engel if and only if  $g_zf_{zz}-f_zg_{zz} \neq 0.$ (We use the notation $\dd{F}{u} = F_u$).

\vspace{2mm}

\begin{lemma}
	Away from $\Sigma$ the characteristic direction is spanned by the vector field
	\[
	C = cZ + eW
	\]
	
	where 
	\[c= g_z (f_{zw} -f_{zy} -f_{zx}+g_zf_y-f_zg_y+f_zf_x)+ f_z (g_{zy}-g_{zw}+g_{zx}-f_zg_x)\]
	\[e= f_zg_{zz}-g_zf_{zz} \]
	
\end{lemma}

\vspace{2mm}

{\sc Proof of Theorem 1,  continued.}

Let us write $S$ for the set of endpoints of singular curves leaving the origin $q_0 = 0 =  (0,0,0,0)$
and lying within a neighborhood for which the normal form holds.  Recall in
all cases that $\Sigma$ is the x-y plane, which is to say the locus $z = w =0$.
Also observe that in all cases that $D$ is transverse to $\Sigma$ which implies
that singular curves, being horizontal, cannot lie inside $\Sigma$.  Since any sufficiently short subarc
of a singular curve remains singular, and outside of $\Sigma$
the only singular curves are the integral curves of $C$,  it follows that the singular curves through $0$ are 
the intergral curves of $C$ which limit to $0$.\\\

{\sc Case 1.}  For our $(2,2,4)$ case  we have $f=z^2$ and $g=zw$, so our characteristic directions are given by the vector field

\[ C = 2z^2w\dd{}{x}+2zw^2\dd{}{y}-2z\dd{}{z}-2w\dd{}{w}\]

For this vector field, $\rho= z^2+w^2$ is a Lyapunov function in a neighborhood of the $xy$-plane. The solution to this vector field passing through $(x(0),y(0),z(0),w(0))$ is 

\[
\dot{w} = -2w; \qquad w = e^{-2t}w(0)
\]
\[
\dot{z} = -2z; \qquad z = e^{-2t}z(0)
\]
\[
\dot{y} = 2zw^2; \qquad y = \frac{-1}{3}e^{-6t}z(0)w^2(0)+y(0)
\]
\[
\dot{x} = 2z^2w; \qquad x = \frac{-1}{3}e^{-6t}z^2(0)w(0)+x(0)
\]	

The flow lines that limit onto  $(0,0,0,0)\in \Sigma$ form the  smooth $2$ dimensional surface $S$, which is
the graph of the map $(z,w) \mapsto (x,y)$ given by  

\[
(x,y) =-\frac{1}{3}\left( wz^2, zw^2 \right)
\]
Thus $S=\{ (x,y,z,w): x=-\frac{1}{3}wz^2, y=zw^2 \}$, the set of singular curves has measure zero.
\\\ \\

{\sc Case 2.}  In the $(2,3,3,4)$ case (A) follows a similar argument taking $f=z$ and $g=z^2w$, the characteristic directions are given now by the vector field

\[ C = -2zw\dd{}{x}- 2z^2w^2\dd{}{y}- 2z\dd{}{z}+2w\dd{}{w} \]

The general solution is given by 
\[
\dot{w} = 2w; \qquad w = e^{2t}w(0)
\]
\[
\dot{z} = -2z; \qquad z = e^{-2t}z(0)
\]
\[
\dot{y} = -2z^2w^2; \qquad y = -2z^2(0)w^2(0)t +y(0)
\]
\[
\dot{x} = -2zw; \qquad x = -2z(0)w(0)t+x(0)
\]

\[
(x,y) = (-zw\ln w ,-z^2w^2\ln w  )
\]

Again  $S = \{  (x,y,z,w): x=-zw\ln w  , y=-z^2w^2\ln w   \}$ is a smooth surface.
\\\ \\

{\sc Case 3.} For the $(2,3,3,4)$ case (B) we have $f=z$ and $g=\frac{1}{3}x_3^3+ x_3x_4^2$, the characteristic directions are given now by the vector field

\[ C = -2z^2\dd{}{x}- 2(\frac{1}{3}z^4 + z^2w^2)\dd{}{y}- 2w\dd{}{z}+2z\dd{}{w}\]

From here we get 
\[
\dot{w} = 2z;
\]
\[
\dot{z} = -2w;
\]

The function  $z^2+w^2$ is constant along solutions.  This implies that
there are no singular curves except the constant curve pass  through the origin:   $S = \{q_0 \}$.
To see this, recall that singular curves cannot lie within $\Sigma$ because $D$ is transverse to $\Sigma$. 
Now, if   a singular curve approaches the origin, then it does so from outside
of $\Sigma$ and eventually enters a neighborhood in which the normal form holds
and there the function  $z^2+w^2$ takes some positive  value $\epsilon$.
From then on, our alleged singular curve must keep this value and so
can never approach the origin which has $z^2 + w^2 = 0$. 

QED \\\ \\

\textbf{Proof of Lemma 1 (following \cite{montgomery2006tour} Chapter 5 )}

We have the characterization of a singular curve of $end$ \ref{s1} - \ref{s3}, we work explicitly condition \ref{s3} where $\dot{\lambda}$ is the tangent vector to the characteristic direction.

Let $\theta_1, \theta_2$ be the local frame for $D^\perp$, so $D^\perp \ni \lambda = \lambda_1\theta_1 + \lambda_2\theta_2 $ in the coordinates $\lambda_1, \lambda_2, x,y,z,w$ ($\lambda_1,\lambda_2$ are fiber coordinates in $D^\perp$), so

\[
d\lambda=d(\lambda_1\theta_1 + \lambda_2\theta_2) = d\lambda_1\theta_1 + d\lambda_2\theta_2 + \lambda_1d \theta_1 +\lambda_2 d \theta_2
\]

since 
\[d\theta_1=(f_x\theta_1 dw + f_y\theta_2 +f_zdz)dw,\] 
\[d\theta_2=(f_g\theta_1 dw + g_y\theta_2 +g_zdz)dw,\] 

we have

\[
d\lambda= d\lambda_1\theta_1 + d\lambda_2\theta_2 + (\lambda_1f_x +\lambda_2g_x)\theta_1dw + (\lambda_1f_y +\lambda_2g_y)\theta_2dw + (\lambda_1f_z +\lambda_2g_z)dzdw
\]

\vspace{2mm}

Note that $\theta_1, \theta_2, dz,dw$ frame all $T^*Q$.  Take the dual frame (so $\theta_i(X_j)=\delta_j^i$  ):

\[X_1=\dd{}{x}, \quad X_2=\dd{}{y}, \quad X_3=\dd{}{z}, \quad X_4=\dd{}{w}-f\dd{}{x}-g\dd{}{y}\]

\vspace{2mm} 

Expand $\dot{\lambda}$ as:

\[
\dot{\lambda} = \dot{\lambda}_1 \dd{}{\lambda_1} + \dot{\lambda}_2 \dd{}{\lambda_2}+ \dot{\gamma}_1X_1 + \dot{\gamma}_2X_2 + \dot{\gamma}_3X_3 + \dot{\gamma}_4X_4
\]

\vspace{2mm}

and evaluate $i_{\dot{\lambda}}d\lambda =0$:

\begin{equation*}
	i_{\dot{\lambda}} d\lambda =\dot{\lambda}_1\theta_1 + \dot{\lambda}_2\theta_2 - \dot{\gamma}_1d\lambda_1-\dot{\gamma}_2 d\lambda_2 + \lambda_1(i_{\dot{\lambda}}d\theta_1) +\lambda_2(i_{\dot{\lambda}}d\theta_2) =0
\end{equation*}

\vspace{2mm}

Since $d \lambda_1, d \lambda_2, \theta_1, \theta_2, dz, dw$ is a basis we have $\dot{\gamma}_1= \dot{\gamma}_2 =0$.

Using $d \theta_1= df\wedge dw$, $d \theta_2= dg\wedge dw $ we get:

\[
i_{\dot{\lambda}}d\theta_1=i_{\dot{\gamma}_3X_3+\dot{\gamma}_4X_4}((f_xdx+f_ydy+f_zdz)\wedge dw)
\]
\[
i_{\dot{\lambda}}d\theta_2=i_{\dot{\gamma}_3X_3+\dot{\gamma}_4X_4}((g_xdx+g_ydy+g_zdz)\wedge dw)
\]

\vspace{2mm}

expressing $dx$ and $dy$ in terms of our frame:

\[
i_{\dot{\lambda}}d\theta_1= -\dot{\gamma}_4f_x\theta_1 -\dot{\gamma}_4f_y\theta_2-\dot{\gamma}_4f_zdz+ \dot{\gamma}_3f_zdw
\]

\[
i_{\dot{\lambda}}d\theta_2= -\dot{\gamma}_4g_x\theta_1 -\dot{\gamma}_4g_y\theta_2-\dot{\gamma}_4g_zdz+ \dot{\gamma}_3g_zdw
\]

So

\begin{align*}
	0=i_{\dot{\lambda}} d\lambda &= \Big( \dot{\lambda}_1-\dot{\gamma}_4(\lambda_1f_x+\lambda_2 g_x) \Big)\theta_1 + \Big( \dot{\lambda}_2-\dot{\gamma}_4(\lambda_1f_y+\lambda_2 g_y) \Big)\theta_2\\
	& \quad + \Big( -\dot{\gamma}_4(\lambda_1f_z+\lambda_2 g_z) \Big)dz + \Big( -\dot{\gamma}_4(\lambda_1f_z+\lambda_2 g_z) \Big)dw
\end{align*}

We can take $(\lambda_1,\lambda_2) = (g_z,-f_z)$ and $\lambda = g_z\theta_1-f_z\theta_2$. Now we have 
\[ d\lambda_1 = g_{zx}(\theta_1 - fdw) + g_{zy}(\theta_2 - gdw) +g_{zz}dz + g_{zw}dw \]
\[ d\lambda_2 = -f_{zx}(\theta_1 - fdw) - f_{zy}(\theta_2 - gdw) -f_{zz}dz - f_{zw}dw \]

and then 

\[ d\lambda_1 \theta_1 = (-g_{zx}f - g_{zy}g + g_{zw})dw\theta_1 + g_{zy}\theta_2\theta_1  +g_{zz}dz\theta_1  \]
\[ d\lambda_2 \theta_2 = -f_{zx}\theta_1\theta_2 -f_{zz}dz\theta_2 + (f_{zx}f + f_{zy}g  - f_{zw})dw\theta_2 \]

In this way:

\[ \lambda d\lambda  =  (g_zk+f_zh)\theta_1\theta_2dw + (g_zf_{zz}-f_zg_{zz})\theta_1\theta_2dz \]

Where $k=f_{zw} -f_{zy} -f_{zx}+g_zf_y-f_zg_y+f_zf_x$ and $h=g_{zy}-g_{zw}+g_{zx}-f_zg_x$

The lemma follows from $\lambda d\lambda = i_C(\theta_1\theta_2dzdw)$

\end{proof}

Summary: 
Morally speaking we are using that $X$ is normally hyperbolic with respect to $\Sigma$. However we have not been able to use normal hiperbolicity to prove Theorem 1 and instead did the work by hand.

We have shown (Theorem 1) that at the points of generic degeneration of an  Engel distribution, the   Sard theorem for the endpoint map  holds, 
and  the set of singular endpoints leaving such a  generic degeneration  point is    either  a $2$  dimensional surface or is that  single  degeneration point. 

What about higher degeneration points? We have left that question open.
\\

\vspace{.5cm}

\bibliography{biblio}{}
\bibliographystyle{plain}

\end{document}